\documentstyle[12pt]{article}
\newtheorem{thm}{Theorem}
\newtheorem{prop}[thm]{Proposition}
\newtheorem{claim}[thm]{Claim}
\newtheorem{problem}{Problem}

\newtheorem{thm-defi}[thm]{Theorem/Definition}
\newtheorem{example}[thm]{Example}
\newtheorem{cor}[thm]{Corollary}
\newtheorem{question}[problem]{Question}
\newtheorem{new-lemma}[thm]{Lemma}

\newtheorem{rem}[thm]{Remark}

\newcommand{\E}{{\cal E}}
\newcommand{\F}{{\cal F}}

\renewcommand{\H}{{\cal H}}

\newcommand{\LB}{{\cal L}}
\newcommand{\R}{{\cal R}}

\newcommand{\M}{{\cal M}}

\newcommand{\Z}{{\cal Z}}

\newcommand{\PP}{{\Bbb P}}

\newcommand{\RealNumbers}{{\Bbb R}}
\newcommand{\Integers}{{\Bbb Z}}
\newcommand{\ComplexNumbers}{{\Bbb C}}
\newcommand{\RationalNumbers}{{\Bbb Q}}

\newcommand{\LongRightArrowOf}[1]{\stackrel{#1}{\longrightarrow}}

\newcommand{\StructureSheaf}[1]{{\cal O}_{#1}}
\newcommand{\EndProof}{\hfill  $\Box$}
\newcommand{\restricted}[2]{#1_{\mid_{#2}}}

\newcommand{\Higgs}{{\cal H}iggs}

\newcommand{\Pic}{{\rm Pic}}
\newcommand{\Sym}{{\rm Sym}}
\newcommand{\Ext}{{\rm Ext}}

\newcommand{\Hom}{{\rm Hom}}

\newcommand{\End}{{\rm End}}

\newcommand{\RelExt}{{\cal E}xt}

\newcommand{\Wedge}[1]{\stackrel{#1}{\wedge}}

\input amssymb.sty

\begin{document}
\begin{center}
\begin{Large}
{\bf 
\noindent
Generators of the cohomology ring of moduli spaces of sheaves on 
symplectic surfaces
}
\end{Large}
\\
Eyal Markman
\footnote{Partially supported by NSF grant number DMS-9802532}
\end{center}


\section{Introduction}
\label{sec-introduction}

Let $\M$ be a moduli space of stable sheaves on a K3 or Abelian surface $S$. 
We express the class of the diagonal in $\M\times \M$ in terms of the
Chern classes of a universal sheaf on $\M\times S$
(Theorem \ref{thm-graph-of-diagonal-in-terms-of-universal-sheaves}). 
Consequently, we obtain generators of the cohomology ring of $\M$. 

In Section \ref{sec-higgs} the surface $S$ is the cotangent bundle of a 
Riemann surface. We recover the result of 
\cite{hausel-thaddeus-I}: {\em The cohomology ring of the
moduli spaces of stable Higgs bundles is generated by the universal classes}.

In Section \ref{sec-stable-cohomology} we concentrate on the case where 
$S$ is a K3 and $\M$ is the Hilbert scheme $S^{[n]}$ of length $n$
subschemes of $S$. 
The construction of the set of generators is uniform in $n$.
Regarding the generators as variables, we get a fixed weighted polynomial 
ring $\R^{[\infty]}$, in infinitely many variables, and 
a natural graded ring homomorphism 
$h: \R^{[\infty]}\rightarrow H^*(S^{[n]},\RationalNumbers)$ onto the 
cohomology ring of $S^{[n]}$. 
The uniform set of generators is minimal in the sense 
that the homomorphism $h$ is injective in degree $\leq n$
(Lemma \ref{lemma-stable-cohomology}). 
In that sense, one may view $\R^{[\infty]}$ as the stable cohomology ring
as $n\rightarrow \infty$. 

Several authors found 
generators for the cohomology ring of a moduli space $\M$ 
of stable sheaves on an algebraic variety $X$.
The K\"{u}nneth factors of the Chern classes of a universal sheaf
are the most natural cohomology classes on $\M$. 
Atiyah and Bott proved that the universal classes generate the cohomology ring 
when $X$ is a curve \cite{AB}. When $X=\PP^2$, 
Ellingsrud and Str{\o}mme \cite{ellingsrud-stromme} 
proved a (stronger) version of Theorem
\ref{thm-graph-of-diagonal-in-terms-of-universal-sheaves} on the level of 
Chow rings. 
Beauville generalized the results of Atiyah-Bott and 
Ellingsrud-Str{\o}mme and proved the analogue of Theorem
\ref{thm-graph-of-diagonal-in-terms-of-universal-sheaves} for
certain moduli spaces, when $X$ is a rational or ruled surface. 
When $X$ is a projective surface, the product $\M\times \M$ is stratified by 
the dimension of the 
extension group $\Ext^2_X(E,F)$ at each point representing a pair 
of sheaves $(E,F)$. Beauville dealt with the case of a trivial stratification.
He observed that, if the extension groups 
$\Ext^2_X(E,F)$ vanish identically, then the diagonal is the
degeneracy locus, of expected dimension, 
of a homomorphism between two vector bundle over $\M\times \M$. 
Theorem \ref{thm-graph-of-diagonal-in-terms-of-universal-sheaves}
follows, in Beauville's case, from Porteous' formula.
In this paper we treat the next simplest case. The diagonal is the
single non-trivial stratum in the stratification of $\M\times \M$,
when the canonical line-bundle of $X$ is trivial. 
Beauville's result was generalized recently (and independently) 
in a different direction: Wei-Ping Li, Zhenbo Qin, and Weiqiang Wang
found generators for the cohomology ring of Hilbert schemes of points
on every smooth projective surface $X$ \cite{lqw}.

The motivation for the results in the present paper came from the study
of the monodromy of moduli spaces of sheaves on a K3 surface $S$
\cite{markman-monodromy}. 
The cohomology ring $H^*(S,\Integers)$ is endowed with a natural symmetric 
pairing ${\displaystyle (\alpha,\beta):=\int_S\alpha^\vee\cup\beta\cup td_S}$,
known as the Mukai pairing. The involution $\alpha\mapsto \alpha^\vee$ acts 
by $-1$ on $H^2(S)$. 
Isometries of $H^*(S,\Integers)$, which preserve the Hodge structure,
lift to auto-equivalences of the derived category of $S$
\cite{orlov}. It is natural to ask if isometries of 
$H^*(S,\Integers)$ lift to symmetries of moduli spaces of stable sheaves on
$S$. Given an algebraic class $v\in H^*(S,\Integers)$, denote by $\M(v)$ the 
moduli space of stable sheaves with Chern character $v$.
In a separate paper \cite{markman-monodromy}, 
we show that the subgroup of the isometry group of $H^*(S,\Integers)$,
stabilizing a Chern character $v$, lifts to a 
subgroup of automorphisms of the cohomology ring
$H^*(\M(v),\Integers)$. Moreover, each such automorphism 
is a monodromy operator for some deformation of
the complex structure of $\M(v)$
(see \cite{markman-reflections} for a special case).
A crucial ingredient, in the proof of
these results, is a formula for the monodromy automorphism 
$\gamma_g$ in terms of the isometry $g$ of $H^*(S,\Integers)$ and 
the universal sheaf on $\M\times S$.
The present paper provides the formula (\ref{eq-class-of-diagonal})
in case $g$ is the identity.


\smallskip
{\em Acknowledgments:} It is a pleasure to acknowledge fruitful conversations 
with D. Huybrechts, S. Mukai, F. Sottile, S. Wang, and K. Yoshioka.
I would like to thank T. Hausel and M. Thaddeus for pointing out
a mistake in an earlier version of our proof of Theorem
\ref{thm-higgs}. I thank K. Ogrady for pointing out a mistake in 
the previous version of Claim \ref{claim-two-calculations-of-m-th-chern-class}.
I thank the referee for correcting a mistake in part
\ref{prop-item-analogue-of-vanishing-of-ext1} of 
Theorem \ref{thm-graph-of-diagonal-in-terms-of-universal-sheaves}
and for several remarks improving the exposition. 

\smallskip
{\em Note:} 
After the completion of this paper, M. Lehn and C. Sorger
computed the ring structure of the cohomology of Hilbert schemes 
for K3 surfaces \cite{LS}.


\section{The diagonal in terms of Chern classes of a universal sheaf}
Let $S$ be a K3 or abelian surface, $\LB$ an ample line bundle on $S$, and
$\M:=\M_\LB(r,c_1,c_2)$ the moduli space of $\LB$-stable
sheaves of rank $r \geq 0$ and Chern classes $c_1$ and $c_2$. 
$\M$ is a smooth and symplectic quasi-projective variety 
\cite{mukai-symplectic-structure}. When the rank is $1$, 
$\M$ is isomorphic to the Hilbert scheme
of zero-dimensional subschemes of length $c_2$. 
If the Chern character $r+c_1+[\frac{(c_1)^2}{2}-c_2]$
is an indivisible class in the cohomology ring $H^*(S,\Integers)$,
then there exists an ample line-bundle $\LB$ for which
$\M_\LB(r,c_1,c_2)$ is complete \cite{yoshioka-chamber-structure}. 
We do not assume completeness of $\M$ in the proof of Theorem
\ref{thm-graph-of-diagonal-in-terms-of-universal-sheaves}. 
Denote by $\pi_{ij}$ the projection from $\M\times S\times \M$
onto the product of the $i$-th and $j$-th factors. Given a flat projective 
morphism $\pi:X\rightarrow Y$ and two coherent sheaves $\E$, $\F$
on $X$, we denote by $\RelExt^i_{\pi}(\E,\F)$ the $i$-th relative extension 
sheaf on $Y$ and by 
\[
\RelExt^!_{\pi}(\E,\F) \ \ = \ \ \sum (-1)^i\RelExt^i_{\pi}(\E,\F) 
\]
the corresponding class in the Grothendieck K-group of $Y$. 

Assume that there exists a universal family over $\M\times S$ 
(this assumption is dropped in Section 
\ref{sec-generalization-to-semi-universal-sheaves}). 

\begin{thm}
\label{thm-graph-of-diagonal-in-terms-of-universal-sheaves}
Let $m$ be the dimension of $\M$ and 
$\E'$, $\E''$  any two universal families of sheaves over 
$\M\times S$. 
\begin{enumerate}
\item
\label{prop-item-class-of-diagonal}
The class of the diagonal in 
$\M\times \M$ is identified by
\begin{equation}
\label{eq-class-of-diagonal}
c_m\left[- \ 
\RelExt^!_{\pi_{13}}\left(
\pi_{12}^*(\E'),\pi_{23}^*(\E'')
\right)
\right].
\end{equation}
\item
\label{prop-item-m-1-class-vanishes}
The class $c_{m-1}\left[- \ 
\RelExt^!_{\pi_{13}}\left(
\pi_{12}^*(\E'),\pi_{23}^*(\E'')
\right)
\right]$ vanishes.
\item
\label{prop-item-analogue-of-vanishing-of-ext1}
The  class 
$c_{m-1}\left[
\RelExt^1_{\pi_{13}}\left(
\pi_{12}^*(\E'),\pi_{23}^*(\E'')
\right)
\right]$ 
vanishes  and $c_{m}\left[
\RelExt^1_{\pi_{13}}\left(
\pi_{12}^*(\E'),\pi_{23}^*(\E'')
\right)
\right]$ is $1-(m-1)!$ times the class of the diagonal.
\end{enumerate}
\end{thm}

An immediate corollary of the Theorem is:
\begin{cor}
\label{cor-kunneth-factors-generate}
If $\M$ is complete, then
the K\"{u}nneth factors of the Chern classes of any universal sheaf $\E$
on $\M\times S$ generate the cohomology ring 
$H^*(\M,\RationalNumbers)$. 
\end{cor}

\noindent
{\bf Proof:}
We use Grothendieck-Riemann-Roch in order to express the class 
(\ref{eq-class-of-diagonal}) as a formula 
of the Chern characters of $\E'$ and $\E''$. 
Given a variety $M$, we denote by 
\begin{eqnarray*}
\ell \ : \ \oplus_{i}H^{2i}(M,\RationalNumbers) & \longrightarrow &
\oplus_{i}H^{2i}(M,\RationalNumbers)
\\
(r+a_1+a_2+\cdots ) & \mapsto & 1+a_1 + (\frac{1}{2}a_1^2-a_2) + \cdots 
\end{eqnarray*}
the universal polynomial map which takes the exponential chern character
of a complex of sheaves to its total chern class. 
Denote by $\pi_i$ the projection from 
$M\times S \times M$ on the $i$-th factor.
Given classes $\alpha_i \in H^*(M\times S,\RationalNumbers)$ 
we set:
\begin{equation}
\label{eq-gamma-delta}
\gamma(\alpha_1,\alpha_2) \ = \ 
c_m\left(
\left\{
\ell\left(
\pi_{13_*}\left[
\pi_{12}^*(\alpha_1)^\vee\cdot\pi_{23}^*(\alpha_2)\cdot
\pi_2^*(td_{S})
\right]
\right)
\right\}^{-1}\right).
\end{equation}
With the above notation, Grothendieck-Riemann-Roch yields the equality
\[
\gamma(ch(\E'),ch(\E'')) \ = \ 
c_m\left\{-\RelExt_{\pi_{13}}^!(\pi_{12}^*(\E'),\pi_{23}^*(\E''))\right\}.
\]
The left hand K\"{u}nneth factors of $\gamma(ch(\E'),ch(\E''))$ are 
contained in the subring generated by the K\"{u}nneth factors of $ch(\E')$.
\EndProof

\begin{rem}
\label{rem-after-prop-identification-of-diagonal}
{\rm
\begin{enumerate}
\item
\label{rem-item-independence-of-choice-of-universal-family}
Note, in particular, that the class (\ref{eq-class-of-diagonal}) is 
independent of the choice of the universal families. We will need
in Section \ref{sec-generalization-to-semi-universal-sheaves} the independence 
of  the class (\ref{eq-gamma-delta}) 
with respect to twists of a universal sheaf by $\RationalNumbers$-Cartier 
divisors. Denote the Neron-Severi group of $\M$ by $NS_{\M}$.
Let us fix $\E$ and consider the map $\gamma_{\E}$, given by 
\[
(L',L'')\ \mapsto \ \gamma(ch(\E)\cdot ch(L'),ch(\E)\cdot ch(L'')),
\]
from the vector space $NS_{\M}\otimes_\Integers \RationalNumbers^{\oplus 2}$ 
to $H^{2m}(\M\times\M,\RationalNumbers)$. It is a polynomial 
map, which is constant along the integral lattice $NS_{\M}\oplus NS_{\M}$.
We conclude that $\gamma_{\E}$ is constant on the whole vector space.
\item
When $m=2$, 
Theorem \ref{thm-graph-of-diagonal-in-terms-of-universal-sheaves}
follows easily from the vanishing of 
$\RelExt^1_{\pi_{13}}\left(
\pi_{12}^*(\E'),\pi_{23}^*(\E'')
\right)$ (\cite{mukai-hodge} Proposition 4.10). 
In general, the sheaf 
$\RelExt^2_{\pi_{13}}\left(\pi_{12}^*(\E'),\pi_{23}^*(\E'')\right)$ 
is a line bundle over the diagonal. 
Hence, the class of the diagonal is represented simply by 
\begin{equation}
\label{eq-class-of-diagonal-is-c-d-of-ext-2}
\frac{(-1)^{m-1}}{(m-1)!}\cdot
c_m\left[\RelExt^2_{\pi_{13}}\left(\pi_{12}^*(\E'),\pi_{23}^*(\E'')
\right)\right]
\end{equation}
(see Example 15.3.1 page 297 in \cite{fulton} concerning the coefficient).
Parts \ref{prop-item-class-of-diagonal} and 
\ref{prop-item-m-1-class-vanishes}
of Theorem \ref{thm-graph-of-diagonal-in-terms-of-universal-sheaves} 
are thus equivalent to the identification of the classes 
$c_i\left[
\RelExt^1_{\pi_{13}}\left(\pi_{12}^*(\E'),\pi_{23}^*(\E'')\right)\right]$, 
$i=m, m-1$ 
in Part \ref{prop-item-analogue-of-vanishing-of-ext1}. 
When $m>2$, it seems easier to prove first 
Parts \ref{prop-item-class-of-diagonal} and \ref{prop-item-m-1-class-vanishes}
of the Theorem. 
\end{enumerate}
}
\end{rem}

\bigskip
Given a complex $E_\bullet$ of coherent sheaves, denote by 
$\Delta_t^{(m+1-t)}(E_\bullet)$
the determinant of the square 
$(m\!+\!1\!-\!t)\times (m\!+\!1\!-\!t)$ matrix 
whose $(i,j)$ entry is the Chern 
class $c_{j-i+t}(E_\bullet)$.  For example, 
\begin{equation}
\label{proteus-determinant}
\Delta_1^{(m)}(E_\bullet) \ := \
\left|\begin{array}{ccccccc}
c_1 & c_2 & c_3 &      &  \cdots &   & c_m\\
1   & c_1 & c_2 &      &         &   & c_{m-1}\\
0   & 1   & c_1 & \cdot&         &   & c_{m-2}\\
    &\cdot&\cdot& \cdot& \cdot&\\
\vdots&   &\cdot     &\cdot      &\cdot         &   & \vdots \\
    &     &     &  \cdot    &  \cdot       & \cdot  & c_2\\
0   &     & \cdots&    & 0       & 1 & c_1
\end{array}\right|(E_\bullet)
\end{equation}
while 
\[
\Delta_m^{(1)}(E_\bullet) \ := \ c_m(E_\bullet).
\]
Note that we have the equality 
$
\Delta_1^{(m)}(E_\bullet) = (-1)^m c_m(-E_\bullet) 
$
(Lemma 14.5.1 and Example 14.4.9 in \cite{fulton}). 

\bigskip
\noindent
{\bf Proof:} (of Theorem 
\ref{thm-graph-of-diagonal-in-terms-of-universal-sheaves}) 
Recall that $\M$ is smooth and symplectic \cite{mukai-symplectic-structure}. 
In particular, it is even dimensional.
First, we construct a complex (\ref{eq-complex-V})
of locally free sheaves on $\M\times \M$, whose $i$-th sheaf cohomology
is $\RelExt^{i+1}_{\pi_{13}}(\pi_{12}^*\E',\pi_{23}^*\E'')$. 
Let $p_\M$ and $p_S$ be the projections from $\M\times S$. 
Set 
\[ 
A_1=\left[p_\M^*p_{\M,*}\left\{\E'\otimes p_S^*\LB^n\right\}
\right]\otimes p_S^*\LB^{-n}.
\]
The evaluation homomorphism from $A_1$ onto $\E'$ 
gives rise to the short exact sequence
\[
0\rightarrow A_0 \rightarrow A_1 \rightarrow \E' \rightarrow 0
\]
for $n$  sufficiently large. 
In particular, $A_1$ is isomorphic to $\oplus^N \LB^{-n}$ on
each fiber.
Below we abuse notation and omit the pullback symbol. 
We have the long exact sequence over $\M\times \M$ 

\begin{eqnarray}
\label{eq-long-exact-seq-of-ext-for-A}
0 & \rightarrow &
\RelExt^0_{\pi_{13}}(\E',\E'') \rightarrow 
\RelExt^0_{\pi_{13}}(A_1,\E'') \rightarrow 
\RelExt^0_{\pi_{13}}(A_0,\E'')  \rightarrow 
\\ \nonumber & \rightarrow  & 
\RelExt^1_{\pi_{13}}(\E',\E'') \rightarrow 
0\hspace{12ex}  \rightarrow 
\RelExt^1_{\pi_{13}}(A_0,\E'')  \rightarrow 
\\ \nonumber & \rightarrow  & 
\RelExt^2_{\pi_{13}}(\E',\E'') \rightarrow 0 \hspace{12ex}
 \rightarrow  \RelExt^2_{\pi_{13}}(A_0,\E'') \rightarrow 0. 
\end{eqnarray}

\noindent
Moreover, $\RelExt^0_{\pi_{13}}(\E',\E'')$ vanishes. 
Denote $\RelExt^0_{\pi_{13}}(A_1,\E'')$ by $V_{-1}$. It is locally free. 
We get the short exact
\begin{equation}
\label{eq-short-exact-ext-seq-E-double-prime}
0 \rightarrow V_{-1} \LongRightArrowOf{\alpha} 
\RelExt^0_{\pi_{13}}(A_0,\E'') \LongRightArrowOf{\beta} 
\RelExt^1_{\pi_{13}}(\E',\E'') \rightarrow 0,
\end{equation}
the isomorphism
\begin{equation}
\label{eq-ext2-is-ext1-with-A0}
\RelExt^1_{\pi_{13}}(A_0,\E'') 
\ \cong \ 
\RelExt^2_{\pi_{13}}(\E',\E''), \ 
\end{equation}
\begin{equation}
\label{eq-ext2-A0-E-vanishes}
\mbox{and the vanishing} \ \ \RelExt^2_{\pi_{13}}(A_0,\E'')  =  0.
\end{equation}

Choose a section $\gamma$ of $\LB^k$ and consider 
the short exact sequence defining the sheaf $Q$ on $S\times \M$
\[
0 \rightarrow \E'' \hookrightarrow \E''\otimes \LB^k \rightarrow Q
\rightarrow 0. 
\]
Set 
\begin{eqnarray*}
V_0 & := & \RelExt^0_{\pi_{13}}(A_0,\E''\otimes\LB^k) \ \ \mbox{and} \\
V_1 & := & \RelExt^0_{\pi_{13}}(A_0,Q). 
\end{eqnarray*}
The sheaf $\RelExt^2_{\pi_{13}}(\E',\E''\otimes\LB^k)$ vanishes for positive
$k$ (use Serre's Duality and stability of $E$ along each fiber $\{E\}\times S$,
$E\in \M$). We claim that the sheaves 
$\RelExt^i_{\pi_{13}}(A_0,\E''\otimes\LB^k)$ vanish for $i\geq 1$
and $k$ sufficiently large. The latter vanishing follows from the analogue of 
the long exact sequence (\ref{eq-long-exact-seq-of-ext-for-A}) 
with $\E''$ replaced by $\E''\otimes\LB^k$.
We get the long exact
\begin{eqnarray}
\label{eq-long-exact-seq-of-ext-for-V}
0 & \rightarrow &
\RelExt^0_{\pi_{13}}(A_0,\E'') \rightarrow 
V_0 \LongRightArrowOf{f}
V_1 \rightarrow 
\\ \nonumber & \rightarrow  & 
\RelExt^1_{\pi_{13}}(A_0,\E'') \rightarrow 0 \longrightarrow 
\RelExt^1_{\pi_{13}}(A_0,Q) \rightarrow 
\\ \nonumber & \rightarrow  & 
0
\end{eqnarray}
(here we used (\ref{eq-ext2-A0-E-vanishes})). 
It follows that $\RelExt^i_{\pi_{13}}(A_0,Q)$ vanishes for $i\geq 1$
and both $V_0$ and $V_1$ are locally free sheaves on $\M\times \M$ 
(the Cohomology and Base Change Theorem). 
Combining (\ref{eq-ext2-is-ext1-with-A0}) and 
(\ref{eq-long-exact-seq-of-ext-for-V}) 
we obtain a description of the diagonal 
as the degeneracy locus of the homomorphism of vector bundles $f$ 
\begin{equation}
\label{eq-f}
0\rightarrow \RelExt^0_{\pi_{13}}(A_0,\E'') \LongRightArrowOf{\iota}
V_0 \LongRightArrowOf{f} V_1 \rightarrow 
\RelExt^2_{\pi_{13}}(\E',\E'')  \rightarrow 0.
\end{equation}

Unfortunately, the homomorphism $f$ is not a regular section
of $\Hom(V_0,V_1)$. The expected codimension of the locus $D_t(f)$, 
where the corank of $f$ is $t$, is given by the formula
\[
\rho(t) = t(r_0-r_1+t). 
\]
In our case, $t=1$ and $r_0-r_1$ is equal to the rank of 
$\RelExt^0_{\pi_{13}}(A_0,\E'')$. 
By (\ref{eq-short-exact-ext-seq-E-double-prime}), 
the rank of the latter is equal to 
\[
(m-2) + \mbox{rank}(V_{-1}),
\]
which is very large. Thus, $D_1(g)$ is empty for a regular section $g$ of 
$\Hom(V_0,V_1)$. 

Combining (\ref{eq-short-exact-ext-seq-E-double-prime})
and (\ref{eq-f}) we get the injective 
composition 
\[
(\iota \circ \alpha) : V_{-1} \ \hookrightarrow \ V_0.
\]
Its image is in the kernel of $f$. 
Set  $g:=\iota \circ \alpha$. We get the complex of vector bundles
\begin{equation}
\label{eq-complex-V}
V_{-1} \ \LongRightArrowOf{g}  \ V_0  \ \LongRightArrowOf{f}  \ V_1
\end{equation}
with sheaf cohomology
$\H_{-1}=0$, $\H_0=\RelExt^1_{\pi_{13}}(\E',\E'')$, and
$\H_1=\RelExt^2_{\pi_{13}}(\E',\E'')$. 
It is easy to check that the sheaf cohomology of the complex
dual to (\ref{eq-complex-V}) is $\H_i=\RelExt^{i+1}_{\pi_{13}}(\E'',\E')$. 
Consequently,
$coker(g^*)=\RelExt^2_{\pi_{13}}(\E'',\E')$ and is hence also supported as
a line bundle on $\Delta$. The Theorem now follows from 
Lemma \ref{lemma-vanishing-of-chern-classes}.
\EndProof


\begin{new-lemma}
\label{lemma-vanishing-of-chern-classes}
Let $V_{-1} \  \LongRightArrowOf{g}  \ V_0  \ \LongRightArrowOf{f}  \ V_1$
be a complex of locally free sheaves of ranks $r_{-1}, r_0, r_1$ 
on a smooth variety $M$ satisfying
\begin{enumerate}
\item
\label{lemma-assumption-sheaf-coho}
The sheaf cohomologies satisfy:
$\H_{-1}:=ker(g)=0$ and $\H_1:=coker(f)$ is supported as a line bundle 
on a smooth subvariety $\Delta$ of pure codimension $m$.

\item
\label{lemma-assumption-alternating-rank}
$m\geq 2$ and $-r_{-1}+r_0-r_1 \ = \ m-2$.
\item
\label{lemma-assumption-coker-g-dual}
$coker(g^*:V_0^*\rightarrow V_{-1}^*)$ 
is also supported as a line bundle 
on $\Delta$. 
\end{enumerate}
Then
\begin{enumerate}
\item
\label{lemma-item-c-m-of-complex-is-Delta}
$c_m(V_\bullet) \ := \ c_m(V_0-V_1-V_{-1}) \ = \ 
\left\{
\begin{array}{ccc}
[\Delta] & \mbox{if} & m \ \mbox{is even}
\\
0 & \mbox{if} & m \ \mbox{is odd}
\end{array}\right.$,
\item
\label{lemma-item-c-m-1-of-complex-vanishes}
$c_{m-1}(V_\bullet) = 0$,
\item
\label{lemma-item-vanishing-of-chern-classes-of-H-0}
The class $c_{m-1}(\H_0)$ vanishes and
$c_m(\H_0)=\left\{
\begin{array}{ccc}
(1-(m-1)!)\cdot [\Delta] & \mbox{if} & m \ \mbox{is even}
\\
(m-1)!\cdot [\Delta] & \mbox{if} & m \ \mbox{is odd}
\end{array}\right.$ 
\end{enumerate}
\end{new-lemma}

\medskip
\noindent
{\bf Proof:} (of Lemma \ref{lemma-vanishing-of-chern-classes}) 
The first degeneracy class of a homomorphism $h:E\rightarrow F$
is given, up to sign, by the Chern class $c_{f-e+1}(F-E)$, where $e$ and $f$ 
are the ranks of the vector bundles $E$ and $F$. 
We begin with a construction, encoding the data of the complex
(\ref{eq-complex-V}) in a homomorphism
(\ref{eq-f-bar}) between two vector bundles, the ranks of which differ by 
$m-1$.

Let $K$ be the kernel of 
\[
(V_{-1}\restricted{)}{\Delta} \longrightarrow (V_{0}\restricted{)}{\Delta}
\]
and denote by $\overline{V}_{-1}$ the quotient
\[
0 \rightarrow K \rightarrow  (V_{-1}\restricted{)}{\Delta}
\rightarrow \overline{V}_{-1}\rightarrow 0.
\]
Let 
\[
b \ : \  X \ \longrightarrow \ M\times M
\]
be the blow up of $\Delta$. Denote by $D\subset X$ the exceptional divisor,
$\iota: D\hookrightarrow X$ its embedding, 
and let $\beta:D\rightarrow \Delta$ be the restriction of $b$.
Define $U$ to be the elementary transform of $b^*V_{-1}$ along $D$
\[
0 \rightarrow U \rightarrow 
b^*(V_{-1})(D) \rightarrow 
(\beta^*\overline{V}_{-1})(D) \rightarrow 0.
\]
$U$ corresponds to the sheaf of meromorphic sections of $b^*V_{-1}$ with, at 
worst, a simple pole along $D$ and ``polar tail'' in $K(D)$. 
We get the short exact sequences
\begin{eqnarray}
\label{eq-U-is-an-extension-of-V-minus-1}
0 \rightarrow b^*V_{-1} \rightarrow 
& U &  
\rightarrow \iota_*(\beta^*K)(D) \rightarrow 0 \ \ \ \ \ \mbox{and}
\\
\label{eq-restriction-of-U-to-D-is-an-extension}
0 \rightarrow \beta^*\overline{V}_{-1} \rightarrow 
& \restricted{U}{D} &
\rightarrow
(\beta^*K)(D) \rightarrow 0
\end{eqnarray}
and the complex (not exact)
\[
U \rightarrow b^*V_0 \LongRightArrowOf{b^*(f)}
b^*V_1.
\]

Proof of part \ref{lemma-item-c-m-of-complex-is-Delta} of Lemma
\ref{lemma-vanishing-of-chern-classes}:
Let 
$\pi \ : \ \PP{U}^* \ \rightarrow \ X$ 
be the bundle of hyperplanes in $U$ and set $u:=r_{-1}-1$ to be the dimension 
of its fibers. 
Denote by $\widetilde{D}$ the inverse image $\pi^{-1}(D)$ in $\PP{U}^*$.
We get the following diagram:
\[
\begin{array}{ccccc}
\PP{U}^* & \LongRightArrowOf{\pi} & X & \LongRightArrowOf{b} & M\times M\\
\hspace{1ex}\ \uparrow \ \tilde{\iota}   & & \hspace{1ex}\ \uparrow \ \iota  
& &\hspace{1ex}\  \uparrow \ e \\
\widetilde{D} & \LongRightArrowOf{p} & D & \LongRightArrowOf{\beta} & \Delta
\end{array}
\]
Let 
\begin{equation}
\label{eq-tautological-bundles-on-PU-dual}
0 \rightarrow \tau \rightarrow \pi^*U \rightarrow q \rightarrow 0
\end{equation}
be the short exact sequence of the tautological sub and quotient bundles.
Note that $q$ is the line bundle 
$\StructureSheaf{\PP{U}^*}(1)$.
Denote by $\overline{V}_0$ the quotient
\[
\overline{V}_0 \ := \ (\pi^*b^*V_0)/\tau
\]
and by $\bar{r}_0$ its rank.
We get the induced homomorphism 
\begin{equation}
\label{eq-f-bar}
\bar{f} \ : \ \overline{V}_0 \ \rightarrow \ \pi^*b^*V_1.
\end{equation}
Note that $\bar{r}_0-r_1=m-1$. 
However, the homomorphism $\bar{f}$ is {\em not} a regular section of 
$\Hom(\overline{V}_0,\pi^*b^*V_1)$. 
We have the following relation in the cohomology ring of $\PP{U}^*$:
\begin{equation}
\label{eq-relation-in-the-cohomology-of-a-projective-bundle}
c_1(q)^{u+1}+\pi^*(c_1(U^*))c_1(q)^{u}+ \cdots + \pi^*c_{u+1}(U^*) \ = \ 0.
\end{equation}
Set
\[
\alpha \ := \ c_1(q)^u - c_1(q)^{u-1}\pi^*(c_1(U)).
\]
Then $\alpha$ satisfies:

i) $\pi_*(\alpha\cdot c_1(q)) \ \ = \ \ 0$ \ \ \ and

ii) The restriction of $\alpha$ to a fiber of $\pi:\PP{U}^*\rightarrow X$ is 
the fundamental class of the fiber.

\noindent
We will prove part 
\ref{lemma-item-c-m-of-complex-is-Delta} of the Lemma by calculating a
cohomology  class in two ways:

\begin{claim}
\label{claim-two-calculations-of-m-th-chern-class}
\begin{enumerate}
\item
\label{claim-item-push-forward-of-porteous-class-is-a-difference}
${\displaystyle 
b_*\pi_*\left[
\alpha\cup\Delta_1^{(m)}(\bar{f})
\right] \ = \ (-1)^m\cdot \left(c_m(V_\bullet)-[\Delta]\right)
}$.
\item
\label{claim-item-push-forward-of-porteous-class-is-0}
${\displaystyle 
b_*\pi_*\left[
\alpha\cup\Delta_1^{(m)}(\bar{f})
\right] \ = \ 
\left\{
\begin{array}{ccc}
0 & \mbox{if} & m \ \mbox{is even}
\\
{[}\Delta{]} & \mbox{if} & m \ \mbox{is odd}
\end{array}
\right.
}$
\end{enumerate}
\end{claim}

\noindent
{\bf Proof:} 
Part \ref{claim-item-push-forward-of-porteous-class-is-a-difference}
of the Claim is proven via a straight-forward calculation. 
We have
\begin{equation}
\label{eq-porteous-class-of-f-bar-is-c-m}
\Delta_1^{(m)}(\bar{f}) \ := \ 
\Delta_1^{(m)}(\pi^*b^*V_1-\overline{V}_0)
\ = \ (-1)^m c_m(\overline{V}_0-\pi^*b^*V_1). 
\end{equation}
The Chern polynomials of the various bundles are related by the following 
equations:
\begin{eqnarray}
\nonumber
c(\tau) & = & \pi^*c(U)/c(q), 
\\
\nonumber
c(\overline{V}_0) & = & \pi^*b^*c(V_0)/c(\tau) \ = \ \pi^*[b^*c(V_0)/c(U)]
\cdot [1+c_1(q)], 
\\
\label{eq-chern-poly-of-bar-V-0-minus-V-one}
c(\overline{V}_0-\pi^*b^*V_1) & = &
\pi^*\left[
\frac{b^*c(V_0)}{b^*c(V_1)\cdot c(U)}
\right]\cdot[1+c_1(q)].
\end{eqnarray}

Combining equations (\ref{eq-porteous-class-of-f-bar-is-c-m}) and 
(\ref{eq-chern-poly-of-bar-V-0-minus-V-one})
yields
\begin{equation}
\label{eq-porteous-class-in-terms-of-classes-on-X}
\Delta_1^{(m)}(\bar{f}) \ = \ 
(-1)^m \left\{
\pi^*c_m[-U+b^*V_0-b^*V_1] \ + \ \pi^*c_{m-1}[-U+b^*V_0-b^*V_1]\cdot c_1(q)
\right\}.
\end{equation}
The choice of $\alpha$ implies the equality
\begin{equation}
\label{eq-integration-along-the-fibers}
\pi_*\left[
\alpha\cup\Delta_1^{(m)}(\bar{f})
\right] \ = \ (-1)^m c_m(-U+b^*V_0-b^*V_1).
\end{equation}
Now, equation (\ref{eq-U-is-an-extension-of-V-minus-1}) yields 
$
c(U) \ = \ b^*c(V_{-1})\cdot c[\iota_*(\beta^*K)(D)]
$
and 
\[
c(-U) \ = \ b^*c(-V_{-1})\cdot c[\iota_*(\beta^*K)(D)]^{-1}.
\]
Equation (\ref{eq-integration-along-the-fibers}) becomes 
\[
\pi_*\left[
\alpha\cup\Delta_1^{(m)}(\bar{f})
\right] \ = \ (-1)^m \cdot 
\left\{b^*c(V_\bullet)\cdot c[\iota_*(\beta^*K)(D)]^{-1}\right\}_{m}
\]
and the projection formula yields
\begin{equation}
\label{eq-push-forward-by-b-composed-with-p-of-alpha-cup-determinantal-locus}
b_*\pi_*\left[
\alpha\cup\Delta_1^{(m)}(\bar{f})
\right] = (-1)^m \cdot \left\{c(V_\bullet)\cdot 
b_*\left(c[\iota_*(\beta^*K)(D)]^{-1}\right)
\right\}_{m}.
\end{equation}
Riemann-Roch without denominators implies  the equality
\[ 
c[\iota_*(\beta^*K)] \ = \  
1+\iota_*P(c_1(\StructureSheaf{D}(D)),\beta^*c_1(K)),
\]
where $P$ is a universal polynomial in two variables
(Theorem 15.3 page 297 in \cite{fulton}). 
The push-forward $b_*\{\iota_*(\beta^*c_1(K)^i)\cdot [D]^j\}$
vanishes for $j<m-1$. If $j=m-1$ and $i>0$, then the push-forward has
degree $> m$. It follows that 
the first $m$ terms of $b_*\{c[\iota_*(\beta^*K)(D)]^{-1}\}$ are equal to 
those of $b_*\{c[\iota_*\StructureSheaf{D}(D)]^{-1}\}$.
Calculating, we get 
\begin{eqnarray}
\nonumber
b_*\{c[\iota_*(\beta^*K)(D)]^{-1}\} & \stackrel{\mbox{up to degree} \ m}{=} & 
b_*\{c[\iota_*\StructureSheaf{D}(D)]^{-1}\} \ = \ 
b_*\{c[\StructureSheaf{X}(D)-\StructureSheaf{X}]^{-1}\}\ = 
\\
\nonumber
b_*\{[1 + D]^{-1}\} & \stackrel{\mbox{up to degree} \ m}{=} & 
b_*\{1-D+D^2 + \cdots \ (-1)^m D^m \} \ = \ 
\\
\label{eq-push-forward-by-b-is-one-minus-class-of-diag}
1 + (-1)^{2m-1}[\Delta] & = & 1-[\Delta].
\end{eqnarray}
Part \ref{claim-item-push-forward-of-porteous-class-is-a-difference} 
of Claim \ref{claim-two-calculations-of-m-th-chern-class}
follows from Equations 
(\ref{eq-push-forward-by-b-composed-with-p-of-alpha-cup-determinantal-locus}) 
and (\ref{eq-push-forward-by-b-is-one-minus-class-of-diag}).


\medskip
Part \ref{claim-item-push-forward-of-porteous-class-is-0}
of Claim \ref{claim-two-calculations-of-m-th-chern-class}
is proven using the Excess Porteous Formula. 
Define $Z$ by the exact sequence
\[
0\rightarrow Z \rightarrow
(b^*V_0\restricted{)}{D} \rightarrow
(b^*V_1\restricted{)}{D} \rightarrow \beta^*\H_1 \rightarrow 0.
\]
$Z$ is a vector bundle over $D$ of rank $r_0-r_1+1$.
Let $\overline{Z}$ be the rank $m$ vector bundle 
\[
\overline{Z} \ := \ \pi^*Z/\tau.
\]
Recall that $\widetilde{D}$ is the inverse image $\pi^{-1}(D)$ in $\PP{U}^*$. 
We get the exact sequence 
\[
0\rightarrow \overline{Z}
\rightarrow  (\overline{V}_0\restricted{)}{\widetilde{D}}
\LongRightArrowOf{\restricted{\bar{f}}{\widetilde{D}}}
\pi^*b^*(V_1\restricted{)}{\Delta} \rightarrow 
p^*\beta^*\H_1 \rightarrow 0.
\]
By the excess Porteous Formula
(\cite{fulton} Example 14.4.7 page 258) 

\begin{equation}
\label{eq-excess-porteous}
\Delta_1^{(m)}(\bar{f}) \ = \
\tilde{\iota}_* c_{\bar{r}_0-r_1}\left[
\overline{Z}^{\vee}\otimes p^*\beta^*\H_1 \ - \ 
\StructureSheaf{\widetilde{D}}(\widetilde{D})
\right]
\ = \
\tilde{\iota}_* c_{m-1}\left[
\overline{Z}^{\vee}\otimes p^*\beta^*\H_1 \ - \ 
\StructureSheaf{\widetilde{D}}(\widetilde{D})
\right]
\end{equation}
as the pushforward to $\PP{U}^*$ of a class in the Chow group 
$A^{m-1}(\widetilde{D})$ of codimension $m-1$ in $\widetilde{D}$. 
{\em Note: Here we use assumption
\ref{lemma-assumption-alternating-rank} 
of Lemma \ref{lemma-vanishing-of-chern-classes}.
}

Given any class $c$ in $A^{m-1}(\widetilde{D})$, we have
\[
b_*\pi_*[\alpha\cdot \tilde{\iota}_*(c)] \ = \ b_*\pi_*[\tilde{\iota}_*
\{\tilde{\iota}^*(\alpha)\cdot c\}] \ = \ 
e_*\beta_*[p_*\{\tilde{\iota}^*(\alpha)\cdot c\}].
\]
Thus, part \ref{claim-item-push-forward-of-porteous-class-is-0}
of the Claim would follow from the equality 
\begin{equation}
\label{eq-the-class-delta-is-the-pushforward-from-D-tilde}
(\beta\circ p)_*\left\{\tilde{\iota}^*(\alpha)\cdot c_{m-1}\left[
\overline{Z}^{\vee}\otimes p^*\beta^*\H_1 \ - \ 
\StructureSheaf{\widetilde{D}}(\widetilde{D})
\right]\right\} \ \ = \ \ 
\left\{
\begin{array}{ccc}
0 & \mbox{if} & m \ \mbox{is even}
\\
{[}\Delta{]} & \mbox{if} & m \ \mbox{is odd}
\end{array}
\right.
\end{equation}

Denote by $D_x$ and $\widetilde{D}_x$
the fibers of $D$ and $\widetilde{D}$ over $x\in \Delta$. 
Let $\delta$ be the restriction of the class 
$\pi^*c_1(\StructureSheaf{D}(D))$ to $A^1(\widetilde{D}_x)$. 
Since $(\beta\circ p):\widetilde{D}\rightarrow \Delta$ 
is a fibration and $\Delta$ is of pure dimension $m$, 
equation (\ref{eq-the-class-delta-is-the-pushforward-from-D-tilde}) 
would follow from the equality of classes on the fiber of $(\beta\circ p)$
\begin{equation}
\label{eq-excess-class-before-integration-along-the-fiber}
\left\{
\tilde{\iota}^*(\alpha)\cdot c_{m-1}\left[
\overline{Z}^{\vee}\otimes p^*\beta^*\H_1 \ - \ 
\StructureSheaf{\widetilde{D}}(\widetilde{D})
\right]
\restricted{\right\}}{\widetilde{D}_x}
\ = \
\left\{
\begin{array}{ccc}
0 & \mbox{if} & m \ \mbox{is even}
\\
{[}\mbox{pt}{]} & \mbox{if} & m \ \mbox{is odd}
\end{array}
\right.
\end{equation}
The exact sequence 
(\ref{eq-restriction-of-U-to-D-is-an-extension}) and
the triviality of $(\beta^*\overline{V}_{-1}\restricted{)}{D_x}$ yields
\[
c(\restricted{U}{D_x}) \ = \ c(\StructureSheaf{D_x}(D_x)) \ = \ 1+\delta.
\]
Denote by $\eta$ the class $c_1(q\restricted{)}{\widetilde{D}_x}$.
From
(\ref{eq-tautological-bundles-on-PU-dual})
we get the equality
\[
c(\restricted{\tau}{\widetilde{D}_x}) \ = \ \frac{1+\delta}{1+\eta}.
\]
Now, the exact sequence
\[
0 \rightarrow \restricted{\tau}{\widetilde{D}_x} \rightarrow 
\pi^*(\restricted{Z}{D_x}) \rightarrow 
\restricted{\overline{Z}}{\widetilde{D}_x} \rightarrow 0,
\]
and the triviality of $\restricted{Z}{D_x}$ imply
\begin{eqnarray*}
c(\restricted{\overline{Z}}{\widetilde{D}_x}) & = &
c(\restricted{\tau}{\widetilde{D}_x})^{-1} \ = \ \frac{1+\eta}{1+\delta}
\ \ \ \mbox{and}
\\
c(\restricted{\overline{Z}^{\vee}}{\widetilde{D}_x}) & = &
\frac{1-\eta}{1-\delta}.
\end{eqnarray*}
Now, 
$c(-\StructureSheaf{\widetilde{D}_x}(\widetilde{D}_x)) = 
\frac{1}{1+\delta}
$. 
Thus, 
\begin{eqnarray*}
\Delta_1^{(m)}(\bar{f}\restricted{)}{\widetilde{D}_x} & = & 
\left\{
c(\restricted{\overline{Z}^{\vee}}{\widetilde{D}_x})\cdot 
c(-\StructureSheaf{\widetilde{D}_x}(\widetilde{D}_x))
\right\}_{m-1} \ = \
\left\{
\frac{1-\eta}{1-\delta^2}
\right\}_{m-1} \ = \
\\
& = & 
\left\{
(1-\eta)(1+\delta^2+\delta^4+\cdots
\right\}_{m-1} \ = \
\\
& = & 
\left\{
\begin{array}{ccc}
-\eta\delta^{m-2} & \mbox{if} & m \ \mbox{is even}
\\
\delta^{m-1} & \mbox{if} & m \ \mbox{is odd.}
\end{array}
\right.
\end{eqnarray*}
The class $\alpha$ restricts to $\widetilde{D}_x$ as $\eta^u-\eta^{u-1}\delta$.
The relation (\ref{eq-relation-in-the-cohomology-of-a-projective-bundle})
yields
\[
-\eta^{u+1}\delta^{m-2} \ = \ 
\eta^u(p^*c_1(U^*\restricted{)}{\widetilde{D}_x})\cdot 
\delta^{m-2} \ = \ -\eta^u\delta^{m-1}.
\]
Hence, 
\[
(\eta^{u}-\eta^{u-1}\delta)\cdot \Delta_1^{(m)}
(\bar{f}\restricted{)}{\widetilde{D}_x} \ = \ 
\left\{
\begin{array}{ccc}
0  & \mbox{if} & m \ \mbox{is even}
\\
\eta^u\delta^{m-1} & \mbox{if} & m \ \mbox{is odd.}
\end{array}
\right.
\]
We have thus proven the equation 
(\ref{eq-excess-class-before-integration-along-the-fiber}) 
of classes on $\widetilde{D}_x$. 
This completes the proof of Claim 
\ref{claim-two-calculations-of-m-th-chern-class}. 
\EndProof

\medskip
Proof of part \ref{lemma-item-c-m-1-of-complex-vanishes}
of Lemma \ref{lemma-vanishing-of-chern-classes}: 
Denote the quotient $b^*V_0/U$ by $E$. It is a vector bundle of rank 
$(m-2)+r_1$. We get the induced homomorphism
\begin{equation}
\label{eq-phi}
\phi \ : \ E \ \longrightarrow \ b^*V_1. 
\end{equation}
Again we calculate a cohomology class in two ways:

\begin{claim}
\label{claim-calculating-m-1-chern-class-in-two-ways}
\begin{enumerate}
\item
\label{claim-item-push-forward-of-porteous-phi-is-chern-class-of-complex}
${\displaystyle
b_*\Delta_1^{(m-1)}(\phi) \ = \ 
(-1)^{m-1}\cdot c_{m-1}(V_\bullet) \ \ \ \mbox{and}
}$
\item
\label{claim-item-push-forward-of-porteous-phi-vanishes}
${\displaystyle
b_*\Delta_1^{(m-1)}(\phi) \ = \ 0.
}$
\end{enumerate}
\end{claim}

\noindent
{\bf Proof:} 
We verify part 
\ref{claim-item-push-forward-of-porteous-phi-is-chern-class-of-complex}
by a straightforward calculation. 
\begin{eqnarray*}
b_*\Delta_1^{(m-1)}(\phi) & := & \Delta_1^{(m-1)}(b^*V_1-E) 
\ = \ (-1)^{m-1}c_{m-1}(E-b^*V_1) \ = \ 
\\
& = &
(-1)^{m-1}b_*\left\{
\frac{b^*c(V_0)}{b^*c(V_1)\cdot c(U)}
\right\}_{m-1} \ = \
\\
& = & 
(-1)^{m-1}b_*\left\{
\frac{b^*c(V_0)}{b^*c(V_1)\cdot [b^*c(V_{-1})\cdot c(\iota_*(\beta^*K)(D))]}
\right\}_{m-1}  \ = \ 
\\
& = & 
(-1)^{m-1}b_*\left\{
b^*c(V_\bullet)\cdot [c(\iota_*(\beta^*K)(D))]^{-1}
\right\}_{m-1} \ = \ 
\\
& = &
(-1)^{m-1}\left\{
c(V_\bullet) \cdot b_*[c(\iota_*(\beta^*K)(D))]^{-1}
\right\}_{m-1} \ = \
\\
& \stackrel{(\ref{eq-push-forward-by-b-is-one-minus-class-of-diag})}{=} &
(-1)^{m-1}\left\{
c(V_\bullet) \cdot (1-[\Delta])
\right\}_{m-1} \ = \ (-1)^{m-1}c_{m-1}(V_\bullet).
\end{eqnarray*}

Part \ref{claim-item-push-forward-of-porteous-phi-vanishes}
 of the claim follows easily from the Excess Porteous Formula.
We have the short exact sequence
\[
0 \rightarrow (\beta^*K)(D) \rightarrow
\restricted{E}{D} \LongRightArrowOf{\restricted{\phi}{D}}
(b^*V_1\restricted{)}{D} \rightarrow  \beta^*\H_1 \rightarrow 0. 
\]
By the Excess Porteous Formula, 
the class $\Delta_1^{(m-1)}(\phi)$ is the push-forward of the class
\[
c_{m-2}\left\{
\left[(\beta^*K)(D)\right]^{\vee}\otimes \beta^*\H_1 - \StructureSheaf{D}(D)
\right\}.
\]
Since $b:D\rightarrow \Delta$ has $(m-1)$-dimensional fibers, 
the homomorphism from the codimension $(m-2)$ Chow group of $D$
to the codimension $(m-1)$ Chow group of $M\times M$ vanishes. 
This completes the proof of Claim 
\ref{claim-calculating-m-1-chern-class-in-two-ways}.
\EndProof

\medskip
Proof of part \ref{lemma-item-vanishing-of-chern-classes-of-H-0} 
of Lemma \ref{lemma-vanishing-of-chern-classes}:
The sheaf $\H_1$ is supported as a line bundle on $\Delta$.
Hence, the first two terms of its total Chern class are
\[
c(\H_1) \ = \ 1 + (-1)^{m-1}(m-1)!\cdot [\Delta] + \cdots
\]
(Example 15.3.1, page 297 in \cite{fulton}). 
We get the equality $c_m(-\H_1) = (-1)^{m}(m-1)!\cdot [\Delta]$. 
If $m$ is even, 
Parts \ref{lemma-item-c-m-of-complex-is-Delta} and 
\ref{lemma-item-c-m-1-of-complex-vanishes} of the Lemma
imply the equalities
\[
[\Delta] \ = \ c_m(V_\bullet) \ = 
\left\{c(-\H_1)\cdot c(\H_0)
\right\}_m \ = \ c_m(-\H_1)+c_m(\H_0) \ = \ 
(m\!-\!1)!\cdot [\Delta] + c_m(\H_0) 
\]
and
\[
0 \ = \ c_{m-1}(V_\bullet) \ = \ 
\left\{c(-\H_1)\cdot c(\H_0)\right\}_{m-1}
\ = \ c_{m-1}(\H_0).
\]
Hence, $c_m(\H_0)=[1-(m-1)!]\cdot [\Delta]$ and $c_{m-1}(\H_0)$ vanishes. 
In case $m$ is odd, the proof is similar. 
This completes the proof of Lemma
\ref{lemma-vanishing-of-chern-classes}.
\EndProof

\section{A generalization for semi-universal sheaves}
\label{sec-generalization-to-semi-universal-sheaves}

We generalize 
Theorem \ref{thm-graph-of-diagonal-in-terms-of-universal-sheaves}
to the case where a universal family $\E$ does not exist globally. 
We define a class $ch(\E)$ replacing the Chern character of a universal 
sheaf (see (\ref{eq-chern-character-of-Q-universal-sheaf})). 
The class of the 
diagonal is then given by the topological formula $\gamma(ch(\E),ch(\E))$,
which is the translation (\ref{eq-gamma-delta})
of (\ref{eq-class-of-diagonal}) via Grothendieck-Riemann-Roch. 


Recall the construction of semi-universal families over $\M\times S$
(\cite{mukai-hodge} Theorem A.5).
We start with a covering $\{U_i\}$ of $\M$, open in the \'{e}tale  
or complex topology, and local universal families
$\E_i$ over $U_i\times S$. We get the direct image vector bundles 
$V_i:=p_*(\E_i\otimes L)$ over $U_i$,  
upon a choice of a sufficiently ample line bundle $L$ on $S$. 
The families $\E_i\otimes V_i^*$ glue to a global semi-universal family
$\F$. Denote by $\rho$ the rank of $V_i$. The projective bundles $\PP{V}_i$
glue to a global projective bundle $p:\PP{V}\rightarrow \M$ 
(which need not be the projectivization of a vector bundle).
The vector bundles $V_i^{\otimes \rho}\otimes (\Wedge{\rho}V_i^*)$ 
glue to a global vector bundle $W$. We get a line bundle 
$\StructureSheaf{\PP{V}}(\rho)$ over $\PP{V}$, which restricts as the
$\rho$-th power of the hyperplane bundle on each $\PP^{\rho-1}$ fiber
and such that $p_*\StructureSheaf{\PP{V}}(\rho)$
is a subbundle of $W$. Denote by $\tilde{\E}$ the 
universal quotient sheaf  of $p^*\F$. $\tilde{\E}$  restricts 
as $\StructureSheaf{}(1)\otimes E$ to 
the $\PP^{\rho-1}\times S$ fiber over $E\in \M$.
Set
\[
\widetilde{V}(1) \ := \ \Hom_{\pi_{1,*}}(p^*\F,\tilde{\E}).
\]
It is a rank $\rho$ vector bundle over $\PP{V}$. We get the relations
\begin{eqnarray*}
p^*\F & \cong & [\widetilde{V}(1)]^*\otimes \tilde{\E} 
\ \ \ \ \ \ \ \ \ \ \mbox{and}
\\
p^*W  & \cong & [\widetilde{V}(1)]^{\otimes\rho}\otimes 
\StructureSheaf{\PP{V}}(-\rho).
\end{eqnarray*}

Denote by $ch(\StructureSheaf{\PP{V}}(1))$ the $\rho$-th root of
$ch(\StructureSheaf{\PP{V}}(\rho))$ with $1$ as the coefficient in
degree zero. We conclude that the class 
$\alpha := ch(\tilde{\E})\cdot ch(\StructureSheaf{\PP{V}}(1))$
is the pullback of a class
\begin{equation}
\label{eq-chern-character-of-Q-universal-sheaf}
ch(\E)
\end{equation}
on $\M\times S$ satisfying 
${\displaystyle
\left[\frac{ch(\F)}{ch(\E)}\right]^\rho = ch(W^*).
}$ 
The Chern classes of $ch(\E)$ have rational coefficients. 
It is easy to check that the class $ch(\E)$ is canonical, up to the product
by the Chern character of a class in $\Pic_\M\otimes \RationalNumbers$. 

Theorem \ref{thm-graph-of-diagonal-in-terms-of-universal-sheaves}
implies that the class
\begin{equation}
\label{eq-the-class-of-the-pull-back-of-diagonal-to-PV-times-PV}
c_m\left\{-\RelExt^!_{\pi_{13}}\left(\tilde{\E},\tilde{\E}\right)\right\}
\end{equation}
is the pullback to $\PP{V}\times \PP{V}$ of the class of the diagonal
in $\M\times \M$. Moreover, the class
(\ref{eq-the-class-of-the-pull-back-of-diagonal-to-PV-times-PV})
is equal to $\gamma(ch(\tilde{\E}),ch(\tilde{\E}))$ as well as to
$\gamma(\alpha,\alpha)$ (see Remark
\ref{rem-after-prop-identification-of-diagonal} part 
\ref{rem-item-independence-of-choice-of-universal-family} 
for the latter equality).
The equality $\alpha=p^*ch(\E)$ implies
the equality
\[
\gamma(\alpha,\alpha) \ = \ (p\times p)^*\gamma(ch(\E),ch(\E)).
\]
We conclude that 
$\gamma(ch(\E),ch(\E))$ is the class of the diagonal in $\M\times \M$ 
(since the pullback $(p\times p)^*$ is an injective homomorphism).

\section{Higgs bundles}
\label{sec-higgs}

Let $\Sigma$ be a smooth compact and connected algebraic curve of genus
$g\geq 2$.
A {\em Higgs bundle} on $\Sigma$ is a pair of a vector bundle $E$
and a $1$-form valued endomorphism $\varphi: E \rightarrow E\otimes K_\Sigma$.
Let $\Higgs$ be the moduli space of rank $r$ 
semi-stable Higgs bundles of degree $d$.
A theorem of Hitchin and Donaldson establishes that the moduli space 
$\Higgs(r,0)$ is homeomorphic to the moduli space of 
semi-simple $GL(r,\ComplexNumbers)$ representations of the
fundamental group of the curve \cite{hitchin-self-duality}. 
The moduli space $\Higgs$ parametrizes
representations of the central extension of $\pi_1(\Sigma)$ by $\Integers$. 
The following is a recent result of T. Hausel and M. Thaddeus:

\begin{thm}
\label{thm-higgs}
\cite{hausel-thaddeus-I}
Assume that $r$ and $d$ are coprime. The cohomology ring 
$H^*(\Higgs,\RationalNumbers)$ is generated by the K\"{u}nneth factors 
of the Chern classes of the universal vector bundle. 
\end{thm}

In a sequel paper, they calculated the relations in the rank $2$ case
\cite{hausel-thaddeus-II}. 
In the rest of this section we deduce Theorem 
\ref{thm-higgs} from Theorem 
\ref{thm-graph-of-diagonal-in-terms-of-universal-sheaves}.
The cotangent bundle $T^*\Sigma$ is a symplectic surface. 
There is a natural bijection between Higgs pairs $(E,\varphi)$ 
on $\Sigma$ and sheaves $F$ on $T^*\Sigma$ with complete support of
pure dimension $1$. $F$ is characterized by the two properties:
1) The support of $F$ is the {\em spectral curve} $C$ of $(E,\varphi)$. 
2) The pushforward $\pi_*(F)$ is isomorphic to $E$. 
If $C$ is smooth, then $F$ is a line bundle on $C$. 
More canonically, the Higgs field $\varphi$ determines a homomorphism from the
sheaf of commutative algebras $\oplus_{i=0}^\infty (K^{-i})$ to $\End(E)$. 
The former is the pushforward of the structure sheaf 
$\StructureSheaf{T^*\Sigma}$ of the cotangent bundle. 
Hence, $E$ is the pushforward of an $\StructureSheaf{T^*\Sigma}$-module $F$.

Let $S:=\PP[T^*\Sigma\oplus \StructureSheaf{\Sigma}]$ be the 
compactification of $T^*\Sigma$. 
Denote the bundle map by $\pi:S\rightarrow \Sigma$ and let 
$D_\infty$ be the section at infinity. 
Choose an ample line bundle $A$ on $S$ and 
consider the moduli space $\M$ of $A$-stable sheaves on $S$ with rank $0$, 
first Chern class 
$c_1:=r\cdot c_1[\pi^*K_\Sigma\otimes \StructureSheaf{S}(D_\infty)]$,
and Euler characteristic $\chi:=d\!+\!r(1\!-\!g)$. 
Stability of the Higgs pair is equivalent to $A$-stability of the
sheaf $F$ \cite{simpson-construction-of-moduli}. 
$\Higgs$ is the Zariski open subset of $\M$ of sheaves with support 
disjoint from $D_\infty$. 
The polarization $A$ can be chosen so that every $A$-semi-stable sheaf is
$A$-stable and $\M$ is thus compact. 
(Indeed,
$A$ can be chosen so that $A\cdot c_1$
and $\chi$ are coprime. Thus, $A$-slope-stability is equivalent to
$A$-slope-semi-stability, where the slope is $\chi/(A\cdot c_1)$ 
and we define slope-stability 
as in \cite{le-potier} Section 2.1).
A universal sheaf $\F_\M$ exists on $\M\times S$ because $r$ and $\chi$
are coprime. (Indeed,
if we set $B_1=\StructureSheaf{S}$ and let $B_2$ be the structure sheaf 
of a fiber of $\pi:S\rightarrow \Sigma$, then $\chi(F\cdot B_1)=\chi$ and
$\chi(F\cdot B_2)=r$, for a sheaf $F$ in $\M$. 
Now apply \cite{huybrechts-lehn} Theorem 4.6.5). 
Denote by $\F$ the restriction of $\F_\M$ to $\Higgs\times S$.

The surface $S$ is not symplectic. Nevertheless, $T^*\Sigma$ is symplectic
and $\F$ is supported on a subscheme of $\Higgs\times T^*\Sigma$. 
It follows that $\F\otimes K_S$ is isomorphic to $\F$. 
Moreover, if $F$ is a stable sheaf in $\Higgs$ and $G$
is a stable sheaf in $\M$, then Serre's Duality yields
\begin{eqnarray}
\label{eq-ext-2-F-G}
\Ext^2_S(F,G) & \cong &
\Hom(G,F\otimes K_S)^* \cong \Hom(G,F)^*, \ \ \ \mbox{and}
\\
\nonumber
\Ext^2_S(G,F) & \cong & 
\Hom(F,G\otimes K_S)^* \cong \Hom(F\otimes K_S^{-1},G)^*  \cong \Hom(F,G)^*.
\end{eqnarray}
Let $\pi_{ij}$ be the projection from $\M \times S \times \Higgs$
onto the product of the $i$-th and $j$-th factors. 
We see that the relative extension sheaves
$\RelExt^2_{\pi_{13}}\left(\pi_{12}^*\F_\M,
\pi_{23}^*\F\right)$ and
$\RelExt^2_{\pi_{13}}\left(\pi_{23}^*\F,\pi_{12}^*\F_\M\right)$
are supported as line-bundles on the diagonal in $\M \times \Higgs$.
The proof of Theorem \ref{thm-graph-of-diagonal-in-terms-of-universal-sheaves}
implies that the class of the diagonal in $\M \times \Higgs$
is given by ${\displaystyle 
c_m\left[-\RelExt^!_{\pi_{13}}\left(\pi_{12}^*\F_\M,
\pi_{23}^*\F\right)\right]}$ (smoothness of $\M$ is not needed because
the diagonal in $\M \times \Higgs$ is smooth and is contained in the smooth
locus of the product. At this stage, the diagonal is
a class in the Chow ring of $\M \times \Higgs$, or in Borel-Moore homology).

Next, we prove that the K\"{u}nneth factors of the Chern classes of $\F$
generate $H^*(\Higgs,\RationalNumbers)$. We will need an additional 
compactification. $\Higgs$ 
admits a natural compactification $\overline{\H}$ with quotient singularities 
(an orbifold) \cite{hausel-compactification,schmitt-compactification}. 
The restriction homomorphism
$H^*(\overline{\H},\RationalNumbers)\rightarrow 
H^*(\Higgs,\RationalNumbers)$ is surjective. 
The surjectivity can be seen using the construction of $\overline{\H}$ 
via symplectic cuts \cite{hausel-compactification}. The natural 
$U(1)$ action on $\Higgs$, which rescales the Higgs field, 
extends to a hamiltonian action on $\overline{\H}$. 
The complement $\Z$ of $\Higgs$ in $\overline{\H}$ is an irreducible 
divisor, which is pointwise fixed under the $U(1)$-action.
The moment map $\mu:\overline{\H}\rightarrow \RealNumbers$, with respect to the
K\"{a}hler symplectic structure of $\overline{\H}$, is a perfect Bott-Morse 
function (\cite{kirwan} Section 5.9). 
Moreover, the value of $\mu$ 
along $\Z$ is the maximal (critical) value \cite{hausel-compactification}. 
Let $y_1, \dots, y_k$ be the critical values of the moment map.
Choose real numbers $x_i$ so that $x_0<y_1<x_1<y_2 < \cdots < y_k<x_k$.
Set $X_i:=\mu^{-1}(x_0,x_i)$. We get the long exact sequences 
\[
\cdots \rightarrow H^*(X_{i+1},X_i) 
\ \rightarrow \ H^*(X_{i+1}) \ \rightarrow \ 
H^*(X_i) \ \rightarrow \ \cdots
\]
For $\mu$ to be perfect means that each of the long exact sequences above 
breaks up into short exact sequences.
Now $X_k=\overline{\H}$ and $X_{k-1}$ is a deformation retract of
$\Higgs$ via the downwards Morse flow \cite{milnor}. 
It follows that the long exact sequence of pairs 
\[
\cdots \rightarrow H^*(\overline{\H},\Higgs) 
\ \rightarrow \ H^*(\overline{\H}) \ \rightarrow \ 
H^*(\Higgs) \ \rightarrow \ \cdots
\]
breaks up into short exact sequences. In particular, we get the required 
surjectivity.

Let $\widetilde{\M}$ be a smooth compactification of $\Higgs$,
which admits morphisms $f_\M$ to $\M$ and $f_{\overline{\H}}$ to
$\overline{\H}$ (which restrict as the identity on $\Higgs$). 
For example, we can choose $\widetilde{\M}$ as a resolution of the closure 
of the diagonal of $\Higgs\times \Higgs$ in $\M \times \overline{\H}$. 
The right hand K\"{u}nneth factors of the diagonal in 
$\widetilde{\M}\times \widetilde{\M}$ restrict to the 
right hand K\"{u}nneth factors of the diagonal in 
$\widetilde{\M}\times \Higgs$.
The K\"{u}nneth factors of the diagonal in 
$\widetilde{\M}\times \widetilde{\M}$ generate $H^*(\widetilde{\M})$.
The restriction homomorphism $H^*(\widetilde{\M})\rightarrow H^*(\Higgs)$ 
is surjective because the restriction from $H^*(\overline{\H})$ factors through
$H^*(\widetilde{\M})$. 
It follows that the K\"{u}nneth factors of the 
diagonal in $\widetilde{\M}\times \Higgs$ generate 
$H^*(\Higgs,\RationalNumbers)$. 
The class of the diagonal
in $\widetilde{\M}\times \Higgs$ is the pullback of the one in
$\M\times \Higgs$. It is given by 
${\displaystyle 
c_m\left[-\RelExt^!_{\pi_{13}}\left(\pi_{12}^*f^*_\M\F_\M,
\pi_{23}^*\F\right)\right]}$. 
Consequently, the K\"{u}nneth factors of the 
Chern classes of $\F$ generate $H^*(\Higgs,\RationalNumbers)$.

It remains to show that the subring of $H^*(\Higgs,\RationalNumbers)$,
generated by the K\"{u}nneth factors of the universal vector bundle 
$\E$, is equal to the subring generated by the K\"{u}nneth factors of $\F$.
The pushforward of $\F$ to $\Higgs\times \Sigma$ is  $\E$.
Denote by $\pi_\H$ and $\pi_S$ the projections from $\Higgs\times S$
and by $p_\H$ and $p_\Sigma$ the projections from $\Higgs\times \Sigma$.
The homomorphism 
\begin{eqnarray}
\label{eq-homomorphism-induced-by-ch-F}
ch(\F) \ \  : \ \ H^*(S,\RationalNumbers) & \rightarrow &
H^*(\Higgs,\RationalNumbers)
\\
\nonumber
\alpha & \mapsto & \pi_{\H,*}[\pi_S^*(\alpha\cdot td_S)\cdot ch(\F)]
\end{eqnarray}
factors through $H^*(T^*\Sigma,\RationalNumbers)$, and hence through 
$H^*(\Sigma,\RationalNumbers)$. 
Given a class $\alpha$ in $H^*(\Sigma,\RationalNumbers)$, we get the equalities
\begin{eqnarray*}
\pi_{\H,*}\left[ch(\F)\cdot \pi_S^*td_S \cdot (1\times \pi)^*p_\Sigma^*(\alpha)
\right]  & = &  
p_{\H,*}\left[(1\times\pi)_*\left\{
ch(\F)\cdot td_{(1\times \pi)} 
\right\}\cdot p_\Sigma^*(\alpha\cdot td_\Sigma)
\right] 
\\
& = & 
p_{\H,*}\left[ch(\E)\cdot p_\Sigma^*(\alpha\cdot td_\Sigma)\right],
\end{eqnarray*}
where the first equality is by the projection formula and
the second by 
Grothendieck-Riemann-Roch. 
It follows that the projection of the image of
(\ref{eq-homomorphism-induced-by-ch-F}) into 
$H^i(\Higgs,\RationalNumbers)$ is equal to the subspace obtained by 
replacing $S$ by $\Sigma$ and $\F$ by $\E$. 
\EndProof

\begin{rem}
{\rm
(Following a suggestion of M. Thaddeus)
Theorem \ref{thm-higgs} was proven in the more general case 
of $L$-valued Higgs bundles, where $L=K_\Sigma(D)$, $D$ an effective divisor
on $\Sigma$ \cite{hausel-thaddeus-I}.
The above argument applies also in the more general case.
When $D>0$, the surface is $S:=\PP[K_\Sigma(D)\oplus \StructureSheaf{\Sigma}]$
and  $K_S\cong \StructureSheaf{S}(-2D_\infty-\pi^*D)$. 
The only difference is that the extension groups in (\ref{eq-ext-2-F-G}) 
vanish. The formula for the diagonal in $\M\times \Higgs$ is obtained,
in this case, using Beauville's result \cite{beauville-diagonal}
(instead of Theorem \ref{thm-graph-of-diagonal-in-terms-of-universal-sheaves}).
}
\end{rem}

\section{The stable cohomology ring}
\label{sec-stable-cohomology}

Let $S$ be a smooth, simply connected, projective surface and 
$\E$ the ideal sheaf of the universal subscheme in $S^{[n]}\times S$.
Denote by $\R^{[n]}$  the weighted polynomial ring generated 
by the vector spaces (of ``variables'')
$H^2(S^{[n]},\RationalNumbers)$ in degree $2$ and $M_{2i}$, $2\leq i\leq n$,
in degree $2i$, where the vector space $M_{2i}$ comes with a fixed isomorphism
with $H^*(S,\RationalNumbers)$. Denote by $\R^{[n]}_k$ the summand of 
weighted degree $k$.
Let $B_{2i}\subset H^*(S^{[n]},\RationalNumbers)$ be the image of
$H^*(S,\RationalNumbers)$ under the homomorphism
\begin{equation}
\label{eq-B-2i}
H^*(S,\RationalNumbers) \ \cong \ H^*(S,\RationalNumbers)^*
\ \LongRightArrowOf{ch(\E)} \ H^*(S^{[n]},\RationalNumbers)
\ \rightarrow \ H^{2i}(S^{[n]},\RationalNumbers).
\end{equation}

We get a natural homomorphism of graded rings
\[
h \ : \ \R^{[n]} \ \rightarrow \ H^*(S^{[n]},\RationalNumbers)
\]
induced by mapping $\R^{[n]}_2$ isomorphically onto $H^2(S^{[n]})$ and
mapping $M_{2i}$ onto $B_{2i}$ via the homomorphism (\ref{eq-B-2i}).

\begin{rem}
\label{rem-invariance-of-chern-character-of-universal-sheaf}
{\rm
Let $S$ be a K3 surface 
and $W$ the reflection group  of $H^2(S^{[n]},\Integers)$ 
with respect to Beauville's bilinear form
\cite{beauville-varieties-with-zero-c-1}. $W$ is a finite index subgroup
of the isometry group of $H^2(S^{[n]},\Integers)$. 
We show in \cite{markman-monodromy} that the Chern character $ch(\E)$
admits a normalization, which is invariant under a diagonal action
of $W$ on $H^*(S^{[n]},\Integers)\otimes H^*(S,\Integers)$. 
The normalization is obtained via multiplication by 1) the square root of
the Todd class of $S$ and 2) the 
Chern character of a class in $\Pic(S^{[n]})\otimes\RationalNumbers$. 
$W$ acts on the cohomology ring of $H^*(S^{[n]},\Integers)$ 
as a subgroup of the monodromy group \cite{markman-monodromy}. 
$W$ acts on $H^*(S,\Integers)$ because the latter contains 
$H^2(S^{[n]},\Integers)$ as a sub-lattice and Beauville's form is the 
restriction of Mukai's pairing \cite{yoshioka-irreducibility}. 
The vector space $H^*(S,\RationalNumbers)$ decomposes as a direct sum of two
irreducible representations of $W$; the trivial character and 
$H^2(S^{[n]},\RationalNumbers)$.
\[ 
H^*(S,\RationalNumbers) \ \ \ \ = \ \ \ \ 
H^2(S^{[n]},\RationalNumbers) \ \oplus \ 1. 
\]
Consequently, $B_{2i}$ is a direct sum of at most two irreducible 
representations of $W$. 
}
\end{rem}

Assume that Corollary \ref{cor-kunneth-factors-generate} holds for
$S$. Then $h$ is surjective and 
G\"{o}ttsche's formula for the Betti numbers of $S^{[n]}$ implies the 
following lower bound on the degree of the relations:

\begin{new-lemma}
\label{lemma-stable-cohomology}
\begin{enumerate}
\item
\label{lemma-item-h-is-surjective}
$h$ is surjective.
\item
\label{lemma-item-h-is-injective-in-degree-leq-n}
$h$ is injective in degree $\leq n = \frac{1}{4}\dim_\RealNumbers(S^{[n]})$.
\item
\label{lemma-item-first-non-trivial-relations}
Let $b_2$ be the second Betti number of $S$. 
The first non-trivial summand of the relation ideal $I\subset \R^{[n]}$ is the 
one-dimensional summand $I_{n+1}$ of degree $n+1$, if $n$ is odd, 
or the $(b_2+2)$-dimensional summand $I_{n+2}$, if $n$ is even. More
generally, the dimension of $I_k$, $k$ even in the range 
$n<k\leq \frac{4n}{3}$, is the sum of the two Betti numbers
\begin{equation}
\label{eq-dimension-of-stable-ideal}
\dim I_k \ \ = \ \ b_{2(k-n-1)}(S^{[n]}) + b_{2(k-n-2)}(S^{[n]}).
\end{equation}
\item
\label{lemma-item-B-2n-vanishes}
If $S$ is a K3, then $H^*(S^{[n]},\RationalNumbers)$ is generated by
$H^2$ and $B_{2i}$, $2\leq i\leq n-1$. 
\end{enumerate}
\end{new-lemma}

In particular, $B_{2i}$ is isomorphic to $H^*(S,\RationalNumbers)$,
for $2i\leq n$, and the relations appear only in degree $> n$. 
Part \ref{lemma-item-first-non-trivial-relations} of Lemma
\ref{lemma-stable-cohomology} seems to suggest that, for $n$ even, $B_{n+2}$ 
is contained in the subring generated by $H^2(S^{[n]})$ and 
$B_{4}$, \dots, $B_{n}$. 
We are led to the following

\begin{question}
Is the ring $H^*(S^{[n]},\RationalNumbers)$
generated by
$H^2$ and the  $B_{2i}$, $2\leq i\leq \lceil\frac{n}{2}\rceil$?
\end{question}

In the K3 case, the answer is affirmative for $n=2$ and $n=3$
(part \ref{lemma-item-B-2n-vanishes} of Lemma \ref{lemma-stable-cohomology}). 
The $n=2$ case follows also from Proposition \ref{prop-verbitzky} 
because the dimension of $H^4(S^{[2]})$ is equal to that of 
$\Sym^2H^2(S^{[2]})$. 
An affirmative answer for $n=4$ seems very plausible 
in view of parts \ref{lemma-item-first-non-trivial-relations}
and \ref{lemma-item-B-2n-vanishes} of Lemma \ref{lemma-stable-cohomology}. 

Lemma \ref{lemma-stable-cohomology} may be interpreted as 
a computation of the stable cohomology ring of the Hilbert schemes. 
Recall that $H^2(S^{[n]},\Integers)$, $n\geq 2$, is the direct sum
\[
H^{2}(S^{[n]},\Integers) \ \ \ \cong \ \ \ 
H^{2}(S,\Integers) \ \oplus \  \Integers\cdot \delta,
\]
where $\delta$ is half the class of the big diagonal. 
Identify the weight $2$ summand of  $\R^{[n]}$, for all $n\geq 2$, 
define the homomorphism $\R^{[n+1]}\rightarrow \R^{[n]}$ by sending 
$M_{2n+2}$ to zero, 
and let $\R^{[\infty]}$ be the inverse limit
\[
\R^{[\infty]} \ := \ \lim_{\infty\leftarrow n}\R^{[n]}. 
\]
$\R^{[\infty]}$ is the weighted polynomial ring generated by
$H^2$ and $M_{2i}$, $i\geq 2$. 
Then the graded kernel ideal of the surjective homomorphism 
$\R^{[\infty]}\rightarrow \R^{[n]} \rightarrow H^*(S^{[n]},\RationalNumbers)$
is trivial in degree $\leq n$.

The proof of Lemma \ref{lemma-stable-cohomology} 
will depend on Lemmas \ref{lemma-dim-R-k} and 
\ref{lemma-stable-betti-number-of-hilbert-scheme}.
The reader may wish to use the table in Section 
\ref{sec-table-of-betti-numbers} to 
test numerically various assertions made in the discussion below.

\begin{new-lemma}
\label{lemma-dim-R-k} 
The dimension of $\R^{[\infty]}_k$ is the coefficient of $t^k$ in the power 
series expansion of the infinite product
\begin{equation}
\label{eq-dim-R-k}
\left[\prod_{m=1}^\infty\frac{1}{(1-t^{2m})^{b_2+1}}\right]
\cdot
\left[
\prod_{m=2}^\infty\frac{1}{1-t^{2m}}
\right].
\end{equation}
\end{new-lemma}

\noindent
{\bf Proof:} 
Let $x_{i,1}$, $1\leq i \leq b_2+1$, be a basis of $\R^{[\infty]}_2$ and let
$x_{i,m}$, $1\leq i \leq b_2+2$, be a basis of $M_{2m}$,
$m\geq 2$. 
Consider the homomorphism
\[
\R^{[\infty]} \ \rightarrow \ \RationalNumbers[t]
\]
sending $x_{i,m}$ to $t^{2m}$. 
Each monomial in the variables $x_{i,m}$ 
appears precisely once in the power series expansion of
\begin{equation}
\label{eq-power-series-in-x}
\left[
\prod_{i=1}^{b_2+1}\prod_{m=1}^\infty\frac{1}{1-x_{i,m}}
\right] 
\cdot
\left[
\prod_{m=2}^\infty\frac{1}{1-x_{b_2+2,m}}
\right].
\end{equation}
Thus, the dimension of $\R^{[\infty]}_k$ is the coefficient of $t^k$
in (\ref{eq-dim-R-k}).
\EndProof

\begin{new-lemma}
\label{lemma-stable-betti-number-of-hilbert-scheme}
The $k$-th Betti number of $S^{[n]}$ is equal to the coefficient of 
$t^k$ in the power series expansion of (\ref{eq-dim-R-k}),
provided $n\geq k$.
\end{new-lemma}

\noindent
{\bf Proof:}
G\"{o}ttsche's formula computes the $k$-th Betti number of $S^{[n]}$
as the coefficient of $t^kq^n$ in the power series expansion of 
\begin{equation}
\label{eq-gottsche-formula}
\prod_{m=1}^\infty
\frac{1}
{(1-t^{2m-2}q^m) \cdot (1-t^{2m}q^m)^{b_2} \cdot (1-t^{2m+2}q^m)}
\end{equation}
\cite{gottsche-formula}.
It will be convenient for us to rewrite it as the product of
$
\frac{1}{1-q}
$
with
\begin{equation}
\label{eq-powwer-serries-in-t-q}
\left[
\prod_{m=1}^\infty\frac{1}{1-t^{2m}q^{m+1}}
\right]
\cdot
\left[
\prod_{m=1}^\infty\frac{1}{(1-t^{2m}q^m)^{b_2}}
\right]
\cdot
\left[
\prod_{m=2}^\infty\frac{1}{1-t^{2m}q^{m-1}}
\right].
\end{equation}

The infinite product (\ref{eq-powwer-serries-in-t-q}) is the image of 
(\ref{eq-power-series-in-x})
under the homomorphism
\begin{equation}
\label{eq-homomorphism-from-R-infinity-to-Q-t-q}
\R^{[\infty]} \ \rightarrow \ \RationalNumbers[t,q]
\end{equation}
given by 
\begin{eqnarray*}
x_{b_2+2,m} & \mapsto & t^{2m}q^{m-1}, \ \ \ \ \ m\geq 2,
\\
x_{b_2+1,m} & \mapsto & t^{2m}q^{m+1}, \ \ \ \ \ m\geq 1,
\\
x_{i,m} & \mapsto & t^{2m}q^m, \ \ \ \ \ 1\leq i \leq b_2, \ \ \mbox{and} \ \ 
m\geq 1.
\end{eqnarray*}
The coefficient of $t^k$ in (\ref{eq-powwer-serries-in-t-q}) is zero, 
if $k$ is odd, and is a polynomial in $q$, 
\begin{equation}
\label{eq-poly-in-q}
c_k(q) \ \ = \ \ 
a_0q^k + a_1q^{k-1}+ \cdots + 
a_{k-\lceil\frac{k}{4}\rceil}q^{\lceil\frac{k}{4}\rceil}
\end{equation}
if $k$ is even. 
The dimension of $\R^{[\infty]}_k$ is the value of $c_k(q)$ at $q=1$. 
Clearly, the degree of $c_k(q)$ is  $k$. 
Thus, the coefficient of $t^kq^n$ in the power series expansion of 
$t^kc_k(q)\cdot\frac{1}{1-q}$ is $c_k(1)$, provided $n\geq k$. 
Hence, the coefficient of $t^kq^n$ in G\"{o}ttsche's formula
is $c_k(1)$, which is the dimension of $\R^{[\infty]}_k$. More generally, 
the difference between the dimension of
$\R^{[\infty]}_k$ and the $k$-th betti number $b_{k}(S^{[n]})$ of $S^{[n]}$
is given by the sum
\begin{equation}
\label{eq-sum-of-leading-coefficients}
\dim\R^{[\infty]}_k-b_{k}(S^{[n]}) \ \ = \ \
\sum_{i=0}^{k-n-1}a_i.
\end{equation}
\EndProof

\medskip
{\bf Proof of Lemma \ref{lemma-stable-cohomology}:} 
\ref{lemma-item-h-is-surjective})
Corollary \ref{cor-kunneth-factors-generate}
implies that $H^*(S^{[n]},\RationalNumbers)$ is generated by 
$H^2(S^{[n]},\RationalNumbers)$ and the $B_{2i}$, $2\leq i\leq 2n$.
We need to show that it suffices to include $B_{2i}$
in the range $2\leq i\leq n$. Clearly, 
$h$ surjects onto the cohomology of degree $\leq 2n$. But 
the cohomology of degree $\leq 2n$
generates the whole cohomology ring (this can be seen, for example, by
Hard Lefschetz). 

\ref{lemma-item-h-is-injective-in-degree-leq-n}) 
We need to prove injectivity in degree $\leq n$.
The dimension of the space $\R^{[n]}_k$, of polynomials of weighted degree
$k$ in $\R^{[n]}$, is equal to the dimension of $\R^{[\infty]}_k$, 
if $k\leq 2n$. 
This dimension is calculated in Lemma
\ref{lemma-dim-R-k} and is shown to be equal to 
the dimension of $H^k(S^{[n]})$ in Lemma
\ref{lemma-stable-betti-number-of-hilbert-scheme}, provided $k\leq n$.

\ref{lemma-item-first-non-trivial-relations})
The dimension of $I_k$ is given by the sum 
(\ref{eq-sum-of-leading-coefficients}).
Let us compute the coefficient $a_i$ in (\ref{eq-poly-in-q}). 
The contribution to the coefficients of $t^kq^n$ in
(\ref{eq-powwer-serries-in-t-q}) comes from monomials in the power
series expansion of (\ref{eq-power-series-in-x}) 
under the homomorphism (\ref{eq-homomorphism-from-R-infinity-to-Q-t-q}).
Let us write these monomials in the form
\[
(x_{b_2+1,1})^d f(\{x_{i,m}\ : \ (i,m)\neq (b_2+1,1)\}),  \ \ \ d\geq 0.
\]
The $(t,q)$ bi-degree of $f(x_{i,m})$ is $(k-2d,n-2d)$ and it satisfies
$k-2d\geq\frac{4}{3}(n-2d)$ (with equality if $f$ is a power of $x_{b_2+1,2}$).
Thus, it suffices to count monomials $f(\{x_{i,m}\ : \ (i,m)\neq (b_2+1,1)\})$ 
whose image under the homomorphism 
(\ref{eq-homomorphism-from-R-infinity-to-Q-t-q})
has bi-degree $(k-n+a,a)$, with $n-a$ even and $a$ in the range 
$0\leq a\leq n$. 
The assumption that $k\leq \frac{4n}{3}$ implies the inequality $a\leq n$. 
The condition that $n-a$ is even is always satisfied, because $k$ is even and 
every monomial is sent by (\ref{eq-homomorphism-from-R-infinity-to-Q-t-q})
to one with even $t$-degree.
Hence, the number of such monomials is the coefficient of 
$z^{k-n}$, under composition of the homomorphism
(\ref{eq-homomorphism-from-R-infinity-to-Q-t-q})
with the homomorphism
\[
\RationalNumbers[t,q] \ \rightarrow \ \RationalNumbers[z,1/z]
\]
sending $t$ to $z$ and $q$ to $1/z$. The image of 
(\ref{eq-powwer-serries-in-t-q}), with its first factor omitted 
(the factor corresponding to $x_{b_2+1,1}$), is the infinite product 
\begin{eqnarray}
\label{eq-poer-series-in-z}
\left[
\prod_{m=1}^\infty\frac{1}{(1-z^m)^{b_2+1}}
\right]
\cdot
\left[
\prod_{m=3}^\infty\frac{1}{1-z^m}
\right]
\ \ = 
\\ 
\nonumber
1+(b_2+1)z+ \frac{(b_2+1)(b_2+4)}{2}z^2 + \ \cdots 
\end{eqnarray}
The infinite product (\ref{eq-poer-series-in-z}) is equal to the one obtained 
from (\ref{eq-dim-R-k}) via the change of variable $z=t^2$ and multiplication 
by $(1-z^2)$. Hence, the coefficient of $z^i$ in (\ref{eq-poer-series-in-z})
is $\dim(\R^{[\infty]}_{2i})-\dim(\R^{[\infty]}_{2i-4})$. 
The dimension of $I_k$, $k$ even in the range 
$n<k\leq \frac{4n}{3}$, is the sum of the coefficients of 
$z^i$, $0\leq i \leq k-n-1$, in the power series expansion of
(\ref{eq-poer-series-in-z}). The equality
(\ref{eq-dimension-of-stable-ideal}) follows.
\bigskip

%
%

\ref{lemma-item-B-2n-vanishes})
Let $Q$ be the quotient of $H^*(S^{[n]},\RationalNumbers)$ by the
subring generated by 
$\oplus_{k=0}^{2n-2} H^k(S^{[n]},\RationalNumbers)$.
Denote by $\overline{B}_{2n}$ the image of $B_{2n}$ in $Q$. 
Assume that the homomorphism (\ref{eq-B-2i}) is normalized as in 
Remark \ref{rem-invariance-of-chern-character-of-universal-sheaf}.
Then $\overline{B}_{2n}$ is a sum of two 
$W$-representations 
$\overline{B}_{2n}=\overline{B}'_{2n}\oplus\overline{B}''_{2n}$.  
$B'_{2n}$ is the image of the trivial $W$-character in the Mukai 
lattice $H^*(S,\Integers)$,
under the normalized homomorphism (\ref{eq-B-2i}) in degree $2n$. 
A priori, $\overline{B}'_{2n}$ is either zero or the trivial 
character and $\overline{B}''_{2n}$ is either zero or isomorphic to 
$H^2(S^{[n]})$. 
The trivial character is spanned by the 
Mukai vector $v=(1,0,1-n)$ of the ideal sheaf of $n$ points.
The Mukai vector $(1,0,1)$ of the trivial line-bundle is not orthogonal to $v$.
Hence, the image of the line spanned by $(1,0,1)$ under (\ref{eq-B-2i}) 
projects {\em onto} $\overline{B}'_{2n}$. 
Similarly, the $W$-orbit of the image of
$(1,0,1)$ in $\overline{B}''_{2n}$ spans $\overline{B}''_{2n}$. 
Since the normalization of (\ref{eq-B-2i}) involves classes in $H^2(S^{[n]})$, 
it suffices to show that the image of $(1,0,1)$ under  (\ref{eq-B-2i})
belongs to the subring generated by $H^k(S^{[n]})$, $k< 2n$. 
Let $\Z\subset S^{[n]}\times S$ be the universal subscheme, $\E$ its ideal 
sheaf, and $p$ the projection to $S^{[n]}$. 
We need to show that the Chern class
$c_n(p_!\E)$ belongs to the subring generated by $H^k(S^{[n]})$, $k< 2n$. 
But $c_n(p_!\E)$ is equal to $c_n(-p_*\StructureSheaf{\Z})$.
The vector bundle $p_*\StructureSheaf{\Z}$ is a direct sum of the trivial 
line-bundle and a rank $n-1$ vector bundle. 
Hence, $c_n(-p_!\E)=0$ and 
$c_n(p_!\E)$ belongs to the subring generated by $H^k(S^{[n]})$, $k< 2n$.
\EndProof

\medskip
Let us recall Verbitsky's result about the subring
$A$ generated by $H^2(S^{[n]},\Integers)$, when $S$ is a K3 surface:

\begin{prop}
\label{prop-verbitzky}
\cite{verbitsky,bogomolov,looijenga-lunts}
Let $X$ be a K\"{a}hler symplectic manifold of dimension
$2n$ and let $A$ be the subalgebra of $H^*(X,\RationalNumbers)$ 
spanned by $H^2(X,\RationalNumbers)$. Then Poincare-Duality restricts
to a perfect pairing on $A$, $H^*(X,\RationalNumbers)=A\oplus A^\perp$, 
and $A$ is  the quotient of 
$\Sym^*H^2(X,\RationalNumbers)$ by the ideal generated by the $n\!+\!1$
powers 
$\alpha^{n+1}$ for all $\alpha\in H^2(X,\RationalNumbers)$ isotropic 
with respect to Beauville's bilinear form.
\end{prop}


\begin{example}
\label{example-n-3}
{\rm
Let us determine, up to one structure constant, the multiplication in 
the cohomology ring of the Hilbert scheme $S^{[3]}$, 
of length $3$ subschemes of a K3 surface (see also \cite{fantechi-gottsche}). 
$H^*(S^{[3]},\RationalNumbers)$ 
is generated by $H^2:=H^2(S^{[3]},\RationalNumbers)$ and $B_4$
(part \ref{lemma-item-B-2n-vanishes} of Lemma \ref{lemma-stable-cohomology}).
The $4$-th betti number of $S^{[3]}$ is $299$. By Verbitzky's result, the 
$276$-dimensional $\Sym^2H^2$ embeds in $H^4(S^{[3]})$ (Proposition 
\ref{prop-verbitzky}). Thus, $H^4(S^{[3]})$ is the direct sum of 
$\Sym^2H^2$ and $B_4''$, where $B_4''$ is 
the 23 dimensional $W$-sub-representation of $B_4$ isomorphic to
$H^2$ (see Remark \ref{rem-invariance-of-chern-character-of-universal-sheaf}). 
Moreover, $B_4'$ is contained in $\Sym^2H^2$. 

The $6$-th betti number of $S^{[3]}$ is $2554$. Again, the $2300$-dimensional
$\Sym^3H^2$ embeds in $H^6(S^{[3]})$. Let us compute the image of 
the homomorphism
\begin{equation}
\label{eq-multiplication-of-H-2-with-B-4}
H^2(S^{[3]})\otimes B_4'' \ \rightarrow \ H^6(S^{[3]}). 
\end{equation}
$H^2(S^{[3]})\otimes B_4''$ decomposes into three irreducible 
$W$-representations: $1\oplus V(2) \oplus \Wedge{2}H^2$, 
where $V(d)\subset \Sym^dH^2$, $d\geq 1$, is the irreducible $W$ 
sub-representation spanned by $d$-th powers of isotropic vectors
(\cite{looijenga-lunts} Proposition 2.14). 
The $W$-decomposition of $\Sym^3H^2$ is $uH^2\oplus V(3)$, 
where $u\in \Sym^2H^2$ is the inverse bilinear form. 
Hence, the image of $H^2(S^{[3]})\otimes B_4''$
in $H^6(S^{[3]})$ intersects the image of $\Sym^3H^2$ trivially. 
The image of the $275$-dimensional 
sub-representation $V(2)$ of $H^2(S^{[3]})\otimes B_4''$
is trivial because $275 > (2554-2300)$. 
Since $H^2$ and $B_4''$ generate the cohomology, the homomorphism 
(\ref{eq-multiplication-of-H-2-with-B-4}) maps both the trivial character and 
the 253-dimensional $\Wedge{2}H^2$ injectively into $H^6(S^{[3]})$. 


$H^6(S^{[3]},\RationalNumbers)$ decomposes as the direct sum of $4$
distinct irreducible representations
\[
H^6(S^{[3]}) \ \ = \ \ uH^2 \oplus V(3) \oplus 1\oplus \Wedge{2}H^2.
\]
Hence, they are pairwise orthogonal with respect to the 
multiplication pairing. The multiplication pairing on $H^6$ 
is determined, up to $4$ constants, 
by the $W$ action. Two of the constants
(for $uH^2 \oplus V(3)$) 
are accounted for by rescaling $H^{12}(S^{[3]})$ and by the structure of 
$A$ in Proposition \ref{prop-verbitzky}. One additional constant 
is eliminated by rescaling $B_4''$. 
We show below that the abstract ring structure is determined 
by the remaining single constant. 

The ring structure is determined by the computation of the 
product of all triples of classes $\alpha_1$, $\alpha_2$, $\alpha_3$,
such that $\alpha_i$ is of pure degree, 
$\deg(\alpha_i)\leq 6$, and $\sum_{i=1}^3\deg(\alpha_i)=12$. 
The case where one of the $\alpha_i$ has degree $6$ has been 
treated. The remaining case is when all three $\alpha_i$ are in $H^4(S^{[3]})$.
If one of the $\alpha_i$ is in the image of $\Sym^2H^2$, we can 
deduce the product from the previous cases. Thus, it remain to determine
the homomorphism $\Sym^3B_4''\rightarrow H^{12}$. But $\Sym^2B_4''$
decomposes as $1\oplus V(2)$ and 
$B''_4$ does not appear in $\Sym^2B_4''$. Hence, 
the homomorphism $\Sym^3B_4''\rightarrow H^{12}$ vanishes.

We have shown that the weighted polynomial ring, generated by $H^2$ and
$B_4''$, surjects onto $H^*(S^{[3]},\RationalNumbers)$ 
and determined its kernel, up to one structure constant. 
\EndProof
}
\end{example}
\section{Appendix: A table of Betti numbers}
\label{sec-table-of-betti-numbers}
The $(k,n)$ entry in the table below is the dimension of $H^k(S^{[n]})$
for a K3 surface $S$. It is calculated via G\"{o}ttsche's formula
(\ref{eq-gottsche-formula}).

\vspace{1ex}
\medskip

\begin{tabular}{|r|r|r|r|r|r|r|r|r|r|r|}  \hline 
\hspace{1ex}
n & 1  & 2   & 3    & 4     & 5      & 6      & 7       & 8        & 9
\\ 
k\hspace{1ex} &    &     &      &       &        &        &         &          &   
\\
\hline
0 & 1  & 1   & 1    & 1     & 1      & 1      & 1       & 1        & 1
\\
2 & 22 & 23  & 23   & 23    & 23     & 23     & 23      & 23       & 23
\\
4 & 1  & 276 & 299  & 300   & 300    & 300    & 300     & 300      & 300
\\
6 &    & 23  & 2554 & 2852  & 2875   & 2876   & 2876    & 2876     & 2876
\\
8 &    & 1   & 299  & 19298 & 22127  & 22426  & 22449   & 22450    & 22450 
\\
10&    &     & 23   & 2852  & 125604 & 147431 & 150283  & 150582   & 150605 
\\
12&    &     & 1    & 300   & 22127  & 727606 & 872162  & 894288   & 897141 
\\
14&    &     &      & 23    & 2875   & 147431 & 3834308 & 4684044  & 4831451
\\
16&    &     &      & 1     & 300    & 22426  & 872162  & 18669447 & 23203208
\\
18&    &     &      &       & 23     & 2876   &  150283 & 4684044  & 84967890
\\ \hline
\end{tabular}

\vspace{1ex}

\medskip
Observe that $b_k(S^{[n]})$ stabilizes for a fixed even $k$ and all $n\geq k$.
For even $n\geq 4$, the difference 
$b_n(S^{[n]})-b_n(S^{[n-2]})$ is $24$, which is the dimension of 
$H^*(S)$ and is also the dimension of $B_n$, by Lemma 
\ref{lemma-stable-cohomology}. 


University of Massachusetts, Amherst, MA 01003 

E-mail: markman@math.umass.edu

\end{document}